\documentclass[11pt]{article}

\usepackage{amsthm}
\usepackage{amsmath}
\usepackage{amssymb}
\usepackage{amsmath, amsthm, bbm}
\usepackage{framed, fullpage}
\usepackage{algorithm}
\usepackage[backref,colorlinks,citecolor=blue,bookmarks=true]{hyperref}

\newcommand{\N}{\mathbb{N}}
\newcommand{\R}{\mathbb{R}}

\newcommand{\ith}{i\text{-th}}
\newcommand{\jth}{j\text{-th}}

\newcommand{\s}[1]{\left\lvert #1 \right\rvert}

\newcommand{\e}{\epsilon}

\newcommand{\sub}{\subseteq}
\newcommand{\sm}{\setminus}

\newtheorem{fact}{Fact}

\newcommand{\eps}{\epsilon}
\newcommand{\eqdef}{\stackrel{\rm def}{=}}
\newtheorem{defn}{Definition}         % A counter for Definition
\newcommand{\BD}{\begin{defn}} \newcommand{\ED}{\end{defn}}
\newcommand{\BE}{\begin{enumerate}} \newcommand{\EE}{\end{enumerate}}
\newcommand{\BI}{\begin{itemize}} \newcommand{\EI}{\end{itemize}}
\newcommand{\calA}{{\cal A}}

\newtheorem{thm}{Theorem}
\newcommand{\BT}{\begin{thm}} \newcommand{\ET}{\end{thm}}
\def\FullBox{\hbox{\vrule width 8pt height 8pt depth 0pt}}
\newcommand{\ourqed}{\;\;\;\FullBox}
\newenvironment{ourproof}{\noindent{\bf Proof:~~}}{\(\ourqed\)}
\newcommand{\BPF}{\begin{ourproof}} \newcommand {\EPF}{\end{ourproof}}
\newenvironment{proofof}[1]{\noindent{\bf Proof of {#1}:~~}}{\(\ourqed\)}
\newcommand{\BPFOF}{\smallskip \begin{proofof}} \newcommand {\EPFOF}{\end{proofof}}
\newcommand{\BEQ}{\begin{equation}} \newcommand{\EEQ}{\end{equation}}
\newcommand{\BEQN}{\begin{eqnarray}}\newcommand{\EEQN}{\end{eqnarray}}
\newtheorem{lem}{Lemma}      % A counter for Lemmas
\newcommand{\BL}{\begin{lem}} \newcommand{\EL}{\end{lem}}
\newtheorem{clm}[lem]{Claim}
\newcommand{\BCM}{\begin{clm}} \newcommand{\ECM}{\end{clm}}
\newtheorem{cor}[thm]{Corollary}      % counter AS FOR Theorems
\newcommand{\BC}{\begin{cor}} \newcommand{\EC}{\end{cor}}
\newcommand{\poly}{{\rm poly}}
\newcommand{\BF}{\begin{fact}} \newcommand{\EF}{\end{fact}}

\date{}

\title{Constructing Near Spanning Trees with Few Local Inspections\footnote{This work is partially based on an extended abstract that appeared in the proceedings of the eighteenth international workshop on randomization and computation
(RANDOM)~\cite{LRR14}.}}
\author{
Reut Levi\thanks{\`{E}cole Normale Sup\`{e}rieure and Universit\`{e} Paris Diderot, France.
This work was done in part when Reut Levi was a PhD student at Tel Aviv University.
  Email: {\tt reuti.levi@gmail.com}.
  Supported in part by ISF grants 246/08, 1147/09, 1536/14, NSF grants  CCF-1217423, CCF-1065125, CCF-1420692 and by ANR grant RDAM}
\and
Guy Moshkovitz\thanks{School of Mathematics, Tel Aviv University, Tel Aviv, Israel 69978.  Email: {\tt guymosko@tau.ac.il}. Supported in part by ISF grant 224/11.}
\and
Dana Ron\thanks{School of Electrical Engineering, Tel Aviv University.
  Tel Aviv 69978, Israel.
  Email: {\tt danar@eng.tau.ac.il}.
  Supported in part by ISF grants 246/08 and 671/13.}
\and
Ronitt Rubinfeld\thanks{CSAIL, MIT.
  Cambridge MA 02139, USA.
School of Electrical Engineering, Tel Aviv University.
  Tel Aviv 69978, Israel.
  Email: {\tt ronitt@csail.mit.edu}.
  Supported in part by NSF grants  CCF-1217423, CCF-1065125, CCF-1420692 and by ISF grant 1536/14.}
  \and
Asaf Shapira\thanks{School of Mathematics, Tel Aviv University, Tel Aviv, Israel 69978. Email: {\tt asafico@tau.ac.il}. Supported in part by ISF Grant 224/11 and a Marie-Curie CIG Grant 303320.}
}

\begin{document}
\maketitle

\begin{abstract}

Constructing a spanning tree of a graph is
one of the most basic tasks in graph theory.
Motivated by
several recent studies of local graph algorithms, we consider the following variant of this problem.
Let $G$ be a connected bounded-degree graph. Given an edge $e$ in $G$ we would like to decide whether $e$ belongs to
a connected subgraph $G'$ consisting of
$(1+\eps)n$ edges (for a prespecified constant $\eps >0$), where the decision for different
edges should be consistent with the same subgraph $G'$. Can this task be performed
by inspecting only
a {\em constant} number of edges in $G$?
Our main results are:

\begin{itemize}

\item We show that if every $t$-vertex subgraph of $G$ has expansion $1/(\log t)^{1+o(1)}$ then one can (deterministically) construct
a sparse spanning subgraph $G'$ of $G$ using few inspections. To this end we analyze a ``local'' version of a famous minimum-weight spanning tree algorithm.

\item We show that the above expansion requirement is sharp even when allowing randomization.
To this end we construct a family of $3$-regular graphs of high girth,
in which every $t$-vertex subgraph has expansion $1/(\log t)^{1-o(1)}$.

\end{itemize}

\end{abstract}

\section{Introduction}

Given a graph $G$, one of the most basic tasks one would like to perform
on $G$ is to find a spanning tree of it or perhaps some other sparse spanning subgraph $G'$.
This task can be easily accomplished using numerous well-known algorithms such as DFS (depth-first search), BFS (breadth-first search)
and more. What all of these algorithms have in common is that in order to
% D 26.1
% tell
decide whether a given edge $e$ belongs to the spanning subgraph $G'$, one has to construct the entire spanning tree. Suppose however that one is not interested
in constructing the entire spanning subgraph $G'$, but rather to be able to ``quickly''
tell  if a given edge $e$ belongs to $G'$ or not.
% By quickly we will mean a {\em constant} number of operations.
By ``quickly'' we mean using a {\em constant\/} number of operations.
% Besides being an interesting problem on its own, our main motivation
% for studying this problem is the recent interest
% in such local algorithms for solving classical problems in graph theory. We elaborate on these related works % in Subsection~\ref{subsec:related}.

Such algorithms are of importance in distributed settings, where processors reside on the
vertices of the graph and the goal is to select as few communication links (edges) as possible
while maintaining connectivity.  Another relevant setting is one in which the graph resides in a centralized database, but different,
uncoordinated, servers have access to it, and are interested in only parts of a common sparse spanning subgraph.
In both cases we would like the decision regarding any given edge to be made after inspecting
only a very small portion of the whole graph, but all decisions must be consistent with the same spanning
subgraph. Such algorithms belong to a growing family of local algorithms
for solving classical problems in graph theory. We elaborate on relevant related works
in Subsection~\ref{subsec:related}.

Let us make a simple observation regarding the task of locally constructing a spanning subgraph.
Note that if one insists on locally constructing a spanning {\em tree} $G'$, then it is easy to see that the
task cannot be performed in general without inspecting
almost all of $G$; that is, this task cannot be achieved using a {\em constant} number of queries to $G$.
%Interestingly, this is in contrast to
%the seemingly related problem of estimating the minimum weight of a
%spanning tree (in an edge-weighted graph) in sublinear-time, which can be
%performed with complexity that does not depend on $n\eqdef
%|V|$~\cite{CRT}.
To see this, observe that if $G$ consists of a single path, then the algorithm must answer positively on all edges, while if $G$ consists of a single cycle then the
algorithm must answer negatively on one edge. However, the
two cases cannot be distinguished without inspecting a
linear number of edges.
% \footnote{It is not hard to extend this argument and show that for any $\epsilon >0$,
% every algorithm for constructing a spanning subgraph with $(1+\epsilon)n$ edges must make
% $\Omega(1/\epsilon)$ queries.}.

So suppose we allow the
algorithm some slackness, and rather than requiring
that $G'$ be a tree, only require that it be relatively sparse,
i.e., contain at most $(1+\eps)n$ edges. Summarizing, the question is then, given $\epsilon > 0$,
for which graphs $G$ can we locally construct a spanning subgraph $G'$ consisting
of $(1+\epsilon)n$ edges, such that given an edge $e \in E(G)$ one can
determine if $e \in G'$ using a constant (that may depend on $\epsilon$ but not on $n$) number
of queries to $G$?

Our main result in this paper, stated informally as Theorem~\ref{theo:informal} below,
shows that the answer to the above question is given by a certain variant of graph expansion, which
we now turn to define. For a graph $G$ and a subset $S \sub V(G)$, we write $\partial_G(S)$
for the set of edges of $G$ with precisely one endpoint in $S$.
We write $\phi_G$ for the (edge) \emph{expansion} of $G$, that is, $\phi_G = \min_S \s{\partial_G(S)}/\s{S}$ where
the minimum is taken over all $S \subseteq V(G)$ of size $1 \leq |S| \leq |V(G)|/2$.
Note that a graph may have small expansion yet contain (large) subgraphs with large expansion.
For example, a vertex-disjoint union of cliques has expansion $0$, yet it contains complete graphs that have the largest expansion possible (for graphs of their order). Let us thus say that
a graph is \emph{$f$-non-expanding} if every $t$-vertex subgraph $H$ satisfies $\phi_H \le f(t)$ (we assume $t > 2 $).
%For example, a path is $f$-non-expanding with $f(x)=1/\floor{x/2}$. On the other extreme, a complete graph is $f$-non-expanding with $f(x)=\ceil{x/2}$.

Our main result in this paper can be informally stated as follows.

\BT[{\bf Informal Statement}]\label{theo:informal}
We have the following dichotomy:
\begin{itemize}

\item If $G$ is $f$-non-expanding for $f(t) \ll 1/\log t$ then one can
locally construct a sparse spanning subgraph of $G$. The algorithm is deterministic.

\item There is a family of $3$-regular graphs $G_n$ that are (roughly) $1/\log t$-non-expanding
so that every (possibly randomized) local algorithm for constructing a sparse spanning subgraph
of $G_n$ must accept every edge of $G_n$.
\end{itemize}
\ET

We refer the reader to Definition~\ref{dfn:SSG-alg} for the precise definition of what it means to locally construct
a sparse spanning subgraph, and to Theorems~\ref{cor:f-non-expand-ub} and~\ref{theo:main} for the precise statements of
the two assertions in Theorem~\ref{theo:informal}.

We note that there are numerous families of graphs that satisfy the condition in the first
item of Theorem~\ref{theo:informal}. For example, it follows from the planar separator theorem of Lipton and Tarjan~\cite{LT79} and its extension by Alon, Seymour and Thomas~\cite{AST90} that planar graphs (and more generally, $H$-minor-free graphs) of bounded degree satisfy the condition of the first item.
Also, observe that since the graphs $G_n$ in the second item of Theorem~\ref{theo:informal} have $3n/2$ edges,
there is no algorithm that can locally construct a spanning subgraph of $G_n$ with $(1+\epsilon)n$ edges
for $\epsilon < 1/2$.

We make two comments regarding the results which appeared in the preliminary conference version of this paper~\cite{LRR14}.
First, it was shown in~\cite{LRR14} that
there are graphs such that any algorithm has to inspect
$\Omega(\sqrt{n})$ edges in order to decide whether a given edge belongs to a spanning subgraph $G'$
containing $(1+\eps)n$ edges, for a constant $\eps$.
However, those graphs resulted from random graphs, which have expansion $\Theta(1)$,
and so could not be used in order to show that the non-expansion requirement
given in the first item of Theorem~\ref{theo:main} cannot be relaxed.
Second, it was shown in~\cite{LRR14} that for certain families of graphs, one can solve
the sparse spanning subgraph problem in time $O(\sqrt{n})$. It is an interesting open problem
to decide whether this can be extended to hold for all bounded-degree graphs.
In fact, it would even be interesting to show that for any bounded-degree graph $G$,
one can find a sparse spanning subgraph using $o(n)$ queries\footnote{Note that if we are allowed to make $\Theta(n)$ queries, then
we can just use the standard BFS or DFS algorithms, which find the entire spanning tree of $G$.}.

%If one is only interested in running time, then
%it is an interesting open problem whether $O(\sqrt{n})$ is also an upper bound for general bounded-degree graphs.\footnote{In~\cite{LRR14} it was shown that %$O(\sqrt{n})$ is an upper bound for a special class
%of graphs, which in particular includes very good expanders.}

% We now turn to discuss some
% further motivation and
% related work.

% D 26.1 Gave motivation earlier on
% \subsection{Motivations and related work}\label{subsec:related}
\subsection{Related work}\label{subsec:related}
As is evident from the above description of the problem,
the model we study here is similar to both classical models, such as distributed and parallel computation,
and to more recent ones. In what follows, we describe these models and some related results,
so as to provide a broad context for our work.

\subsubsection{Local algorithms for other graph problems}\label{subsec:local}
The model of {\em local computation algorithms} as considered in
this work, was defined by Rubinfeld et al.~\cite{RTVX} (see
also Alon et al.~\cite{ARVX12}). Such algorithms for maximal
independent set, hypergraph coloring, $k$-CNF and maximum
matching are given in~\cite{RTVX,ARVX12,MRVX12,MV13}.
This model generalizes other models that have been studied
in various contexts, including locally decodable codes
(e.g.,~\cite{STV99}), local decompression~\cite{DLRR13},
and local filters/reconstructors
\cite{ACC+08,SS10,B08,KPS08,JR11,CS06a}. Local computation
algorithms that give approximate solutions for various
optimization problems on graphs, including vertex cover,
maximal matching, and other packing and covering problems,
can also be derived from sublinear time algorithms for
parameter estimation~\cite{PR,MR,NO,HKNO09,YYI}.

The model of local computation is related to several other models, including property testing and online algorithms.
To give a notable example, Mansour et al.~\cite{MRVX12} proposed a general scheme for converting a large family of online algorithms into local computation algorithms, consequently, improving the complexity of hypergraph 2-coloring and $k$-CNF in the local computation model.

In the related field of local reconstructors, Campagna et al.~\cite{CGR13} study the property of connectivity.
Namely, under the promise that the input graph is almost
connected, their reconstructor provides oracle access to
the adjacency matrix of a connected graph which is close to
the input graph. We emphasize that our model is
different from theirs, in that they allow the addition of
new edges to the graph, whereas our algorithms must provide
spanning graphs whose edges are present in the original
input graph.

\subsubsection{Distributed and parallel algorithms}\label{subsec:dist}
The name {\em local algorithms} is also used in the
distributed context~\cite{MNS95,NS95,Lin92}. As observed by
Parnas and Ron~\cite{PR}, local distributed algorithms can
be used to obtain local computation algorithms as defined
in this work, by simply emulating the distributed algorithm
on a sufficiently large subgraph of the graph $G$. However,
while the main complexity measure in the distributed
setting is the number of rounds (where it is usually
assumed that each message is of length $O(\log n)$), our
main complexity measure is the number of queries performed
on the graph $G$. By this standard reduction, the bound on
the number of queries (and hence running time) depends on
the size of the queried subgraph and may grow exponentially
with the number of rounds. Therefore, this reduction gives
meaningful results only when the number of rounds is
significantly smaller than the diameter of the graph.

While the problem of computing a spanning graph has not been studied in the distributed model, the problem of computing a minimum-weight spanning tree is a central one in this model.
Kutten and Peleg~\cite{KP98}
provided an algorithm that works in $O(\sqrt{n} \log^* n +
D)$ rounds, where $D$ denotes the diameter of the graph.
Their result is nearly optimal in terms of the complexity
in $n$, as shown by Peleg and Rubinovich~\cite{PR00} who
provided a lower bound of $\Omega(\sqrt{n}/\log n)$ rounds
(when the length of the messages must be bounded).

Another problem studied in the distributed setting that is
related to the one studied in this paper, is finding a
sparse spanner. The requirement for spanners is much
stronger since the distortion of the distance should be as
small as possible. Thus, to achieve this property, it is
usually the case that the number of edges of the spanner is
super-linear in $n$. Pettie~\cite{Pet10} was the first to
provide a distributed algorithm for finding a low
distortion spanner with $O(n)$ edges without requiring
messages of unbounded length or $O(D)$ rounds. The number
of rounds of his algorithm is $\log^{1+o(1)}n$. Hence, the
standard reduction of~\cite{PR} yields a local algorithm
with a trivial linear bound on the query complexity.

%\subsubsection{Parallel algorithms}\label{subsec:para}
The problems of computing a spanning tree and a minimum
weight spanning tree were studied extensively in the
parallel computing model as well (see, e.g.,~\cite{BC05}, and the
references therein). However, these parallel algorithms
have time complexity which is at least logarithmic in $n$
and therefore do not yield an efficient algorithm in the
local computation model. See~\cite{RTVX,ARVX12} for further
discussion on the relationship between the ability to
construct local computation algorithms and the parallel
complexity of a problem.

\subsubsection{Local cluster algorithms}\label{subsec:clust}
Local algorithms for graph theoretic problems have also
been given for PageRank computations on the web
graph~\cite{JW03,Ber06,SBC+06,ACL06,ABC+08}. Local graph
partitioning algorithms have been presented
in~\cite{ST04,ACL06,AP09,ZLM13,OZ13}, which find subsets of
vertices whose internal connections are significantly
richer than their external connections in time that depends
on the size of the cluster that they output.
For instance, Andersen and Peres~\cite{AP09} provide an algorithm which, given a starting vertex $v$, finds a cluster of $v$ of small conductance,
whose complexity depends on the volume of the cluster it outputs but has only polylogarithmic dependence in the size of the graph.
However, even when the
size of the cluster is guaranteed to be small, it is not
obvious how to use these algorithms in the local
computation setting where the cluster decompositions must
be consistent among queries to all vertices.

\subsubsection{Other related sublinear-time approximation algorithms for graphs}\label{subsec:sublin}
The problem of estimating the weight of a minimum-weight
spanning tree in sublinear time was considered by Chazelle,
Rubinfeld and Trevisan~\cite{CRT}. They describe an
algorithm whose running time depends on the approximation
parameter, the average degree and the range of the weights,
but does not directly depend on the number of vertices.

\subsection{Organization}

The rest of the paper is organized as follows. In Section~\ref{sec:prel}
we formally define the local sparse spanning subgraph problem which we consider in this paper.
In Section~\ref{sec:upper} we prove the first item of Theorem~\ref{theo:informal}, which is formally stated as Theorem~\ref{cor:f-non-expand-ub}. The proof of this part has two main steps. In the first one, we show that if $G$ is $f$-non-expanding with $f \ll 1/\log t$ then one can remove from $G$ only
% $o(n)$
a relatively small number of
edges and thus partition it into connected components of size $O(1)$ each. We then show that if a graph can be so partitioned, then one can solve on it the local spanning subgraph problem by executing a ``localized'' version of Kruskal's~\cite{Kruskal56} famous algorithm for finding minimum-weight spanning tress\footnote{Recall that if $G$ is a graph with weights assigned to its edges, then Kruskal's algorithm finds a spanning tree of minimal {\em total} weight}.

The proof of the second paper of Theorem~\ref{theo:informal}, which is the more challenging part of this paper, is
given in Section~\ref{sec:lower} and formally stated as Theorem~\ref{theo:main}. It establishes that the $1/\log t$-non-expansion
requirement from the first item of Theorem~\ref{theo:informal} is essentially tight. What we show is that there
are graphs which are (about) $1/\log t$-non-expanding, and have the property that any local algorithm for constructing
a spanning subgraph using a constant number of queries must accept {\em every} edge of the graph. To prove this result we describe a construction
of certain extremal graphs that might be of independent interest. These are $3$-regular graphs, that on one
hand have unbounded \emph{girth}\footnote{As usual, the girth of a graph is the minimum length of a cycle in it.}, but on the other hand are (about) $1/\log t$-non-expanding.
% D 26.1: moved to an earlier point in the paper
% We note that in a
% preliminary version of this paper~\cite{LRR14}, it was shown that there are graphs for which one cannot
% locally construct a sparse
% spanning subgraphs. However, those graphs resulted from random graphs, and thus had expansion $\Theta(1)$.
% Theorem~\ref{theo:main} improves upon this by constructing such hard graphs that are
% $1/\log t$-non-expanding.

We make no serious attempt to optimize the constants obtained in the various statements. In fact, the $f$-non-expansion requirements in our upper and lower bound results (Theorems~\ref{cor:f-non-expand-ub} and~\ref{theo:main}), which are about $(1/\log t)(1/\log\log t)^2$ and $(1/\log t)(\log\log t)^2$ respectively, can each be improved by replacing the $(\log\log t)^2$ term by $(\log\log t)^{1+o(1)}$.
We opted for proving our results with the slightly weaker bounds in order to simplify the presentation.
We henceforth write $\log(\cdot)$ for $\log_2(\cdot)$.

\section{Preliminaries}\label{sec:prel}

Let us now give the precise definition of the algorithmic problem we are addressing in this paper.
As in most cases where one tries to design a local/distributed/sublinear algorithm, we will assume
that the input graph is given via an {\em oracle access to its incidence-list representation}, meaning the following:
First, we assume that the input graph $G=(V,E)$
is given via incidence-lists representation, that is, for each vertex $v \in V(G)$, there is an
ordered list of its neighbors in $G$. Second, the algorithm is supplied
with integers $n$ and $d$, that represent the number of vertices, and an upper bound on the degrees of vertices of $G$.
Finally, given a pair $(v,i)$ with $1 \leq v \leq n$ and $1 \leq i \leq d$, the oracle either returns the
$i^{th}$ neighbour of $v$ (in the incidence list representation) or an indication that $v$ has less than $i$ neighbours.
We will assume that each vertex $v$ has an id, $id(v)$, where there is a full order over the ids. We will think of the ids of vertices
in the graphs simply as the integers $\{1,\ldots,n\}$.
We now turn to formally define the algorithmic problem we consider in this paper.

\renewcommand{\labelenumi}{\theenumi}
\renewcommand{\theenumi}{\roman{enumi}.}

\BD %[Local algorithm for sparse spanning graphs]
\label{dfn:SSG-alg}
An algorithm $\calA$ is an {\em $(\e,q)$-local sparse spanning graph algorithm} if, given $n,d \ge 1$ and oracle access to the incidence-lists representation of a connected graph $G=(V,E)$ on $n$ vertices and degree at most $d$,
%parameters
%$n\geq 1$ and $0 \leq \delta < 1$ and given
%query access to the incidence-lists representation of a
%connected graph $G=(V,E)$,
%the algorithm $\calA$ provides
it provides
query access to a subgraph $G'=(V, E')$ of $G$ such that:
%the following hold:
%\BI
\BE
\item $G'$ is connected. % (with probability $1$.
\item $\s{E'} < (1+\eps)\cdot n$ with probability at least
    $2/3$ (over the
    internal randomness of $\mathcal{A}$).
\item $E'$ is determined by $G$ and the internal randomness
    of $\mathcal{A}$.
\item $\calA$ makes at most $q$ queries to $G$.
\EE
% Namely, on
By ``providing query access to $G'$'' we mean that
on input $(u, v)\in E$,
    $\calA$ returns whether $(u, v) \in E'$ and for any
    sequence of queries, $\calA$ answers consistently with
    the same $G'$.

An algorithm $\calA$ is an {\em $(\e,q)$-local sparse spanning graph
algorithm for a family of graphs $\mathcal{C}$}
if the above conditions hold, provided that the input graph
$G$ belongs to $\mathcal{C}$.
\ED

\renewcommand{\labelenumi}{\theenumi}
\renewcommand{\theenumi}{\arabic{enumi}.}

We note that the choice of the required success probability being $2/3$ is of course arbitrary
and can be replaced by any probability smaller than $1$. Having said this, let us stress that the positive results
we obtain here (i.e., the algorithmic results) in Theorem~\ref{cor:f-non-expand-ub} are deterministic
(i.e., hold with probability $1$), whereas our lower bound in Theorem~\ref{theo:main} holds for \emph{any} positive success probability.
We also note that even though Definition~\ref{dfn:SSG-alg} considers only  the number of queries performed by the algorithm, our algorithm in  Theorem~\ref{cor:f-non-expand-ub} runs in time polynomial in the number of queries, and in particular, independent of $n$.

We are interested in
local algorithms that have query complexity which is
independent of $n$, namely, that perform a constant number
of queries to the graph (for each edge they are queried on)
and whose running time (per queried edge) is small as well.
In the next section, we show that the family of graphs that are
$f$-non-expanding with $f \ll 1/\log t$ have a local sparse spanning graph
algorithm. In the following section, we will show that one cannot prove
such a result when $f$ is only slightly larger.

%As stated in the next theorem, we cannot obtain a constant
%query complexity for general bounded degree graphs.

%\BT[\cite{LRR14}]\label{theo:LB} Any local sparse spanning
%graph algorithm has query complexity $\Omega(\sqrt{n})$.
%This result holds for graphs with a constant degree bound
%$d$ and for constant $0 \leq \eps \leq 2d/3$ and $0 \leq
%\delta < 1/3$. \ET

% TODO: are multigraph relevant to us in this paper?
%\footnote{Graphs are allowed to have self-loops and multiple edges, but for our problem we may assume that there are no self-loops and multiple-edges (since  theanswer on a self-loop can always be negative, and the same is true for all but at most one among a set of parallel edges).}
% If the graph is edge-weighted, then the weight of
% the edge is returned as well.

\section{Upper bound}\label{sec:upper}

In this section we prove the following theorem, which formalizes the first assertion of Theorem~\ref{theo:informal}.

\BT\label{cor:f-non-expand-ub}
For every $C$ there is a function $q:\R_+\times\N\to\N$ so that for every $\epsilon >0$ there is an $(\e,q(\e,d))$-local
sparse spanning graph algorithm for the family of $f$-non-expanding graphs with
\begin{equation}\label{eq:f}
f(x) = \frac{C}{\log x\cdot (\log\log x)^2} \;,
\end{equation}
where $d$ is the input degree-bound. Furthermore, the algorithm is deterministic. %for the input graph.
\ET

%Note that for a constant $\e$, the algorithm in Theorem~\ref{cor:f-non-expand-ub} (Algorithm~\ref{alg:local-kruskal} below) has query complexity $q$ and %running time that are both polynomial in $d$. %$d^{O(1)}$.

\subsection{Decomposition of non-expanding graphs}

The first step in the proof of Theorem~\ref{cor:f-non-expand-ub} is a decomposition result stated in Lemma~\ref{lemma:CSS} below.
It shows that if $G$ is $f$-non-expanding, with $f$ as in Equation~(\ref{eq:f}),
then $G$ can be decomposed into connected components of bounded size by removing only $\e n$ edges.
This extends a result of~\cite{CSS09} that applies for somewhat larger $f$.
As mentioned earlier, there are many families of graphs which are $f$-non expanding with $f$ as in Equation~(\ref{eq:f}).
For example, planar graphs of bounded degree are $f$-non-expanding with $f=O(1/\sqrt{x})$ by the famous planar separator theorem of Lipton and Tarjan~\cite{LT79}. More generally, a result of Alon, Seymour and Thomas~\cite{AST90} implies that
for any fixed $H$, the family of $H$-minor-free graphs of bounded degree is $f$-non-expanding with $f=O(1/\sqrt{x})$.
Hence, Lemma~\ref{lemma:CSS} applies to these families of graphs in particular.
%This means, that in such cases we can also set $k=\poly(1/\epsilon)$ in step $(i)$ of Algorithm~\ref{alg:local-kruskal}, thus obtaining
%a much faster running time.
We note that the reason why the bound in Lemma~\ref{lemma:CSS} is doubly exponential in $\epsilon$ is that we insist on assuming that $f$ is very close to the threshold of $1/\log x$ (which by Theorem~\ref{theo:main} is essentially tight).
For example, the details of the proof of Lemma~\ref{lemma:CSS} show that if $f=x^{-c}$ for some $0 < c < 0$ (as is the case
with planar graphs, say), then the bound can be improved to polynomial in $1/\epsilon$.
We note that in such cases we can also set $k=\poly(1/\epsilon)$ in step $1$ of our algorithm (Algorithm~\ref{alg:local-kruskal} below), thus obtaining a much
% faster running time.
more efficient algorithm.

\BL\label{lemma:CSS}
If $G$ is an $n$-vertex $f$-non-expanding graph with
$f(x) = C/\log x(\log\log x)^2$, then one can remove $\e n$ edges from $G$
so that each connected component of the remaining graph is of size at most
$2^{2^{2(C/\e) + 3}}$.
\EL

\BPF
%\BPFOF{Lemma~\ref{lemma:CSS}}
First, we claim that any $f$-non-expanding $n$-vertex graph $G=(V,E)$ has
a subset $S \subset V(G)$ of size $n/3 \le \s{S} \le (2/3)n$
and expansion $\phi_{G}(S) \eqdef \s{\partial_G(S)}/\s{S} \le f(n/3)$.
For the proof we iteratively construct subsets $S_1,\ldots,S_k \sub V(G)$ as follows.
To obtain $S_i$, we consider the induced subgraph $G_i = G[V \sm \bigcup_{j=1}^{i-1} S_j]$ %G-\bigcup_{j=1}^{i-1} S_j$
and let $S_i \sub V(G_i)$ satisfy $\s{S_i} \le n_i/2$ and $\phi_{G_i}(S_i) \le f(n_i)$, where $n_i = \s{V(G_i)}$.
We stop once $S \eqdef \bigcup_{i=1}^k S_i$ is of size $\s{S} \ge n/3$.
Note that %$n_k \le 2n/3$ so $\s{S} \le n/3 + n_k/2 \le 2n/3$.
$$\s{S} \le \sum_{i=1}^{k-1} \s{S_i} + n_k/2 = (n+\sum_{i=1}^{k-1} \s{S_i})/2 \le 2n/3 \;.$$
It remains to bound $\phi_{G}(S)$.
Observe that every edge in the edge boundary $\partial_G(S)$ is a member of some edge boundary $\partial_{G_i}(S_i)$. Hence,
$$\frac{\s{\partial_G(S)}}{\s{S}} \le
\frac{\sum_{i=1}^k \s{\partial_{G_i}(S_i)}}{\s{S}} =
\sum_{i=1}^k \frac{\s{S_i}}{\s{S}} \phi_{G_i}(S_i)
\le \max_{1\le i \le k} \phi_{G_i}(S_i) \le \max_{1\le i \le k} f(n_i) = f(n_k) \le f(n/3) \;,$$
where in the last inequality we used the fact that $n_k \ge n-\s{S} \ge n/3$.
This proves our claim.

Fix an integer $k \geq 50$ and let $G=(V,E)$ be any $f$-non-expanding graph on $n \geq k/3$ vertices.
Consider the following process; take any subset $S \subset V$ of size $n/3 \le \s{S} \le n/2$ and expansion $\phi_{G}(S) \le 2f(n/3)$
(whose existence follows from the claim above), remove all its outgoing edges and proceed recursively on the two induced subgraphs $G[S]$ and $G[V \sm S]$, which are $f$-non-expanding as well.
% We do nothing
The recursion stops
whenever we reach a graph on at most $k$ vertices.
It is clear that at the end of this process, the edges removed from $G$ leave a graph whose connected components have at most $k$ vertices each. Let $r_k(G)$ be the number of edges
removed by the above process. We will shortly prove that if $G$ has $n$ vertices, then $r_k(G) \leq Cn/\ln\ln(k/3)$. Hence, setting $k=2^{2^{2(C/\e) + 3}} \ge \max\{50, 3 \cdot e^{e^{C/\e}}\}$ %50\cdot e^{e^{C/\e}}$
enables us to remove no more than $\epsilon n$ edges and break $G$ into connected components of size at most $k$, thus
proving the lemma.

In order to facilitate an inductive proof, it will be more convenient to prove the following slightly stronger claim:
\begin{equation}\label{eq:ind-claim}
r_k(G) \leq \beta(n) \eqdef \frac{Cn}{\ln\ln(k/3)} - \frac{Cn}{\ln\ln n}\;.
\end{equation}
% Put
Set $h(x) = x/\ln\ln x$ and $f^*(x) = f(x)/C=(\log x)^{-1}(\log \log x)^{-2}$.
First, we %need
establish some properties of $h$.
It is easy to verify that the derivative of $h$ is
$h'(x) = (\ln\ln x)^{-1} - (\ln x)^{-1}(\ln\ln x)^{-2}$,
and moreover, $h''(x) \le 0$ for $x \ge 20$.
It follows that for every $n \ge 50$ and $n/3 \le s \le n/2$ we have
\begin{equation}\label{eq:h-bound}
h(n) - h(n-s) \le s \cdot h'(n-s)  \le h(s)-s \cdot 2f^*(n/3)
\end{equation}
where in the first inequality we used the concavity of $h$ on the interval $(20,\infty)$,
and in the second inequality we used the fact that $n \geq 50$ and $n/3 \leq s \leq n/2$ and that in this range
$$
(\log(n/3)(\log\log(n/3))^2/2 \geq \ln(2n/3)(\ln\ln(2n/3))^2\geq \ln(n-s)(\ln\ln(n-s))^2\;.
$$

We prove Equation~(\ref{eq:ind-claim}) by induction on $n$. For the base case(s) where $(k/3 \le)$ $n \le k$ we have that $\beta(n) \ge \beta(k/3) = 0 = r_k(G)$, as needed.
For the induction step we have
\begin{align*}
r_k(G) &\le \max_{\substack{S \subset V:\\ n/3 \le \s{S} \le n/2}} \s{S} \cdot 2f(n/3) + r_k(G[S]) + r_k(G[V \sm S]) \\
&\le \max_{n/3 \le s \le n/2} s \cdot 2f(n/3) + \beta(s) + \beta(n-s)\\
&= C\Big(n/\ln\ln (k/3) + \max_{n/3 \le s \le n/2} s \cdot 2f^*(n/3) - h(s) - h(n-s) \Big)\\
&\le C\Big(n/\ln\ln (k/3) - h(n) \Big) = \beta(n)
\end{align*}
where the first inequality follows from the definition of the process described in the second paragraph of the proof,
the second inequality follows from the induction hypothesis since $k/3 \le s,n-s \le n-1$, and in the third inequality we used~(\ref{eq:h-bound}) since $n \ge k \ge 50$.
This completes the proof of Equation~(\ref{eq:ind-claim}).
%\EPFOF
\EPF

\subsection{The algorithm}

The algorithm we design in order to prove Theorem~\ref{cor:f-non-expand-ub} is based on Kruskal's minimum-weight spanning tree algorithm~\cite{Kruskal56}. The idea is to assign weights to the edges of the graph in a way that will determine some fixed spanning tree $T$. The algorithm will always accept the edges of $T$ but will also accept a few other edges.
We will pick the weights of the edges in a way that will make it possible to determine the edges of a sparse spanning subgraph in a ``local'' fashion.

Recall that Kruskal's algorithm for finding a minimum-weight spanning tree in a
weighted connected graph works as follows. First it sorts the edges of the graph $e_1,\dots e_m$ from minimum to maximum weight (breaking ties arbitrarily).
%Let this order be $e_1,\dots e_m$.
It then goes over the edges in this order, and adds $e_i$ to the spanning tree if and only if it does not close a cycle with the previously selected edges.
Put differently:

\BF\label{fact:Kruskal}
Edge $e$ is picked by Kruskal's algorithm if and only if for any cycle $C$ of $G$ containing $e$, the edge $e$ does not
have the largest weight among the edges of $C$.
\EF

It is well known (and easy to verify)
that if the weights of the edges are distinct, then there is a single minimum-weight spanning tree in the graph. For an unweighted graph $G$, consider the order
defined over its edges by the order of the ids of the vertices. Namely, we define
a ranking $r$  of the edges as follows:
$r(u,v) < r(u',v')$ if and only if $\min\{id(u),id(v)\} < \min\{id(u'),id(v')\}$
or $\min\{id(u),id(v)\} = \min\{id(u'),id(v')\}$
and $\max\{id(u),id(v)\} < \max\{id(u'),id(v')\}$.
If we run Kruskal's algorithm using the rank $r$ as the weight function (where
there is a single ordering of the edges), then
we obtain a (well-defined) spanning tree of $G$.

While the local algorithm described next (Algorithm~\ref{alg:local-kruskal})
is based on the aforementioned
global algorithm, it does not exactly emulate it, but rather emulates a certain {\em relaxed\/}
version of it which can be executed {\em locally}. In particular, it will answer YES for every edge selected by the global algorithm
(ensuring connectivity), but may answer YES also on edges not selected by the global algorithm.
We will thus need to show that it does not answer YES on too many edges that are not selected
by the global algorithm.

In the description and analysis of the algorithm we will use
the following standard notation; for a vertex
$v \in V$ and an integer $k$, we denote by $C_k(v,G)$ %(or simply $C_k(v)$ if $G$ is clear from context)
the subgraph of $G$ induced by the set of vertices at distance at most $k$ from $v$.

%In the description and analysis of the algorithm we will use
%the following standard notation: we denote the distance between
%two vertices $u$ and $v$ in a graph $G$ by $d_G(u,v)$. For vertex
%$v \in V$ and an integer $k$, let $\Gamma_k(v,G)$ denote
%the set of vertices at distance at most $k$ from $v$ and
%let $C_k(v,G)$ denote the subgraph of $G$ induced by
%$\Gamma_k(v,G)$. When the graph $G$ is clear from the
%context, we shall use the shorthand $d(u,v)$, $\Gamma_k(v)$
%and $C_k(v)$ for $d_G(u,v)$, $\Gamma_k(v,G)$ and
%$C_k(v,G)$, respectively.

\medskip

\begin{algorithm}
\caption{{\bf (Kruskal-based Algorithm)}}
(The algorithm works for some fixed $\epsilon>0$.)

\textbf{Input:} $n,d \ge 1$, query access to a graph $G$ on $n$ vertices and degree at most $d$.

\textbf{Queried edge:} $(x,y) \in E(G)$.
\BE
\item Set $k=2^{2^{2(C/\e) + 3}}$.
\item Perform a BFS to depth $k$ from $x$, thus obtaining the subgraph $C_k(x,G)$.
% induced by $\Gamma_k(x)$ in $G$.
\item If $(x,y)$ is the edge with largest rank on some cycle in $C_k(x,G)$ then answer NO;\\ Otherwise, answer YES.
\EE
\label{alg:local-kruskal}
\end{algorithm}
%%
% As stated in the next theorem, the complexity of our algorithm
% does not depend on the size of the graph.
%We next show:
%\BT\label{thm:kruskal-hyper}\label{theo:UB}
%Algorithm~\ref{alg:local-kruskal}, given $\e \ge 0$,
%when run with $k = \rho(\epsilon)$,
%is a local sparse spanning graph algorithm
% for the family of $f$-non-expan graphs with a degree bounded by $d$.
% The query complexity and running time of Algorithm~\ref{alg:local-kruskal}
% are  $O(d^{\rho(\epsilon)+1})$, and its success probability is $1$.
%\ET

%\BPF  %OF{Theorem~\ref{thm:kruskal-hyper}}

\BPFOF{Theorem~\ref{cor:f-non-expand-ub}}
We will show that if $G=(V,E)$ is $C/\log x(\log \log x)^2$-non-expanding then Algorithm~\ref{alg:local-kruskal} is an $(\epsilon,q(\epsilon,d))$-local sparse spanning subgraph algorithm, where $q(\epsilon,d)=d^{k+1}$ with $k$ being the constant from step $1$ of the algorithm.
By the description of Algorithm~\ref{alg:local-kruskal} it directly follows that
the algorithm is deterministic and that its answers are consistent with a connected subgraph $G'$.
Indeed, if $T$ is the tree returned by Kruskal's algorithm, then Fact~\ref{fact:Kruskal} and step 3 of Algorithm~\ref{alg:local-kruskal} guarantee that each edge of $T$ will be accepted by Algorithm~\ref{alg:local-kruskal}. Observe that the number of queries to $G$
performed by Algorithm~\ref{alg:local-kruskal} is at most $d^{k+1}$.
We now complete the proof by showing that the algorithm returns YES on fewer than $(1+\epsilon)n$ edges.

Let $R$ be a set of at most $\epsilon n$ edges whose removal disconnects $G$ into connected components of size at most $k$. The existence of such a set is guaranteed by Lemma~\ref{lemma:CSS}.
Let $G_R$ be the graph obtained by removing $R$ from $G$; that is, $G_R=(V,E \setminus R)$.
We note (crucially) that while the analysis of the algorithm uses properties of $G_R$, the algorithm does not actually compute $R$. We will now show that $G'$ does not contain a cycle of $G_R$. Since $|R| \leq \epsilon n$, this means
that $G'$ has fewer than $(1+\epsilon)n$ edges.

Let $\sigma$ be a cycle in $G_R$.
%and suppose $\sigma$ belongs to a connected component of $G_R$.
%and consider the connected component of $G_R$ which contains $\sigma$.
Suppose $(w,v)$ is the edge of $\sigma$ with the largest rank. Since the connected components of $G_R$ are of size at most $k$, we infer that $\sigma$ has at most $k$ vertices, implying that $C_k(w,G)$ contains $\sigma$. It follows that on query $(w,v)$ the algorithm will return NO. Thus, $G'$ does not contain $\sigma$.
%
%For each vertex $u$ in $G$, let $\widetilde{C}(u) = (\widetilde{V}(u), \widetilde{E}(u))$ denote the component of
%$u$ in $G_R$. We claim that for each $u$, the graph $G'$ does not contain a cycle on the subgraph induced on $\widetilde{V}(u)$,
%implying that $G'$ contains no cycle of $G_R$. Since $|\widetilde{V}(u)| \le k$, the diameter of $\widetilde{C}(u)$ is at most $k-1$.
%This implies that $C_k(u)$ contains $\widetilde{C}(u)$.
%Let $\sigma$ be a cycle in $\widetilde{C}(u)$ and let $e=(w,v)$ be the edge in $\sigma$ with the largest rank.
%Since $\widetilde{C}(u) = \widetilde{C}(v) = \widetilde{C}(w)$ it follows that on query $(w, v)$ the algorithm returns NO.
\EPFOF
%\medskip
%
%\BPFOF{Theorem~\ref{cor:f-non-expand-ub}}
%Immediate from Theorem~\ref{theo:UB} using the $k$ given by Theorem~\ref{theo:CSS}.. This implies that the number of queries (as well as the running time) is at most some $q(\e,d) = d^{2^{2^{O(1 + C/\e)}} + 1}$.
%\EPFOF

\section{Lower bound}\label{sec:lower}

% Roman numerals in enumerate environment
\renewcommand{\labelenumi}{\theenumi}
\renewcommand{\theenumi}{\roman{enumi}.}

%\mnote{D: Changed introductory text}
%By Theorem~\ref{theo:LB}, there are graphs for which every local sparse spanning graph algorithm
%must perform a number of queries that grows with $n$. However, the graphs used in the
%corresponding lower-bound construction contain
%%
% It is not hard to see that Theorem~\ref{theo:LB} implies that if $G$ consists of a random $d$-regular % graph on $t_0$ vertices connected by a path to the rest of the $n$ vertices, then with high probability %any local sparse spanning graph algorithm for $G$
% it requires $\sqrt{t_0}$ queries.
% Taking $t_0=t_0(n)$ to be $\omega(1)$, we have that $G$ requires an unbounded number of queries.
% Unfortunately, $G$ also has
%large subgraphs with expansion $\Omega(1)$.
%%
% The following result essentially closes the gap between the \emph{hereditary} expansion of the upper bound and lower bound graphs.
 The next theorem shows that there are graphs for which any local sparse spanning graph
 algorithm must perform a number of queries that grows with $n$,
       yet these graphs are $f$-non-expanding with $f(x)$ only slightly larger than $1/\log x$.
        This is essentially the best one can hope for in light of Theorem~\ref{cor:f-non-expand-ub}.

\BT\label{theo:main}
For infinitely many $n$, there is an %connected $3$-regular
$f$-non-expanding $n$-vertex graph $G$ with
$$f(x) = \frac{1}{\log x}\cdot (70\log\log x)^2$$
%$\epsilon \leq \frac14$,
such that every $(\frac12,q)$-local sparse spanning graph algorithm
for the graphs isomorphic to $G$ satisfies $q \ge \log\log(n)/8000$.
\ET

\subsection{A regular non-expanding graph}

The main result in this subsection (stated in
Lemma~\ref{lemma:cube-cycle}) is a construction of
regular non-expanding graphs that we will use in
Subsection~\ref{sec:LB-proof} to prove
Theorem~\ref{theo:main}. A main ingredient is a result
from~\cite{MS15} showing that, roughly speaking, there are
graphs % of unbounded degree
that simultaneously have large
girth and small hereditary expansion (in fact, small edge
separators). While the degree of these graphs may grow with $n$, their maximum
degree is at most a constant times their average degree.
We will use this in order to construct a regular graph
with similar properties. We note that the regularity
condition is crucial for proving Theorem~\ref{theo:main}.
The following theorem was proved in~\cite{MS15}.

\BT[\cite{MS15}]\label{theo:cube}
For any $n,k$ with $2 \le k \le \frac{1}{648}\log\log n$ there is an $n$-vertex graph
$G=G_{n,k}$ satisfying:
\begin{enumerate}
\item $G$ has average degree at least $k$ and maximum degree at most $6k$.
\item $G$ has girth at least $\log n/(6k)^2$.
\item For every $t$-vertex subgraph $H$ of $G$ that is not a forest, there exists a  subset $S \sub V(H)$ of size $(1/3)t \le \s{S} \le (2/3)t$ such that
$$\s{\partial_H(S)} \le  \frac{t}{\log t}\cdot (\log\log t)^2 \;.$$
%\item Every $t$-vertex subgraph $H$ of $G$
%that is not a forest
%with $t\ge \log n/(6k)^2$
%satisfies
%\begin{equation}%\label{eq:expansion}
%$$\phi_H \le \frac{1}{\log t}\cdot 6(\log\log t)^2 \;.$$
%\end{equation}
\end{enumerate}
\ET

We note that each of the parameters in
Theorem~\ref{theo:main} is quantitatively essentially
optimal (see~\cite{MS15} for further discussion).

%The main result in this section is a construction of $3$-regular graphs that satisfy properties similar to those in Theorem~\ref{theo:cube}.
The main result in this section is the following.

\BL\label{lemma:cube-cycle} For any $n_0$ there is a
connected graph $G_\circ$ on $n \ge n_0$ vertices
satisfying:
\begin{enumerate}
\item $G_\circ$ is $3$-regular. %and $n = 2\pmod{4}$.
\item $G_\circ$ has girth at least $\log\log(n)/2000$.
\item $G_\circ$ is $f$-non-expanding with $f(x) = (1/\log x)\cdot (4\log\log x)^2$.
\end{enumerate}
\EL

For the proof we will need the weighted version of the
well-known vertex separator theorem for trees. For
completeness, we give a short proof below.
\BCM\label{clm:tree} Let $T=(V,E)$ be a tree, and let
$w:V\to\R^+$ be a nonnegative weight function over the vertices of $T$.
There is a
vertex $v\in V$ whose removal disconnects $T$ into
connected components of weight at most $w(V)/2$
each.\footnote{For a subset $X \sub V$ we write
$w(X)=\sum_{v\in X} w(v)$.}
\ECM

%\begin{proof}
\BPF
Start a walk in $T$ from an arbitrary vertex, in each step moving from a vertex $u$ to a neighbor $u'$ if the weight of the tree rooted at $u'$, when the edge $(u,u')$ is removed, is strictly greater than $w(V)/2$.
Since $T$ has no cycles and since the walk never reverts the last step taken, the walk eventually stops at some vertex $v$. This means that when $v$ is removed from $T$, the weight of the tree rooted at each of the neighbors of $v$ is at most $w(V)/2$. Since these trees are the connected components resulting
from the removal of $v$, we are done.
\EPF
% \end{proof}

% \begin{proof}[Proof of Lemma~\ref{lemma:cube-cycle}]
\BPFOF{Lemma~\ref{lemma:cube-cycle}}
Set $k=\log\log m/648$ and let $G_{m,k}$ be the graph from Theorem~\ref{theo:cube}, where we take $m$ to be large enough such that
%$k \ge 2n_0$ $(\ge 2)$.
$k \ge \min\{n_0, 2\}$.
%(notice that, by Theorem~\ref{theo:cube}, we necessarily have $2 \le (\log\log m)/648$).
We note that in the rest of the proof we will use the inequality
\begin{equation}\label{eq:relations}
%6k \le m^{1/5} \qquad\text{ and }\qquad \log m/(6k)^2 \ge (6k)^2 \;,
(6k)^4 \le \log m
\end{equation}
which holds since $m$ is sufficiently large. %(automatically from the restriction on $k$). %can be easily verified (notice .
As is well known, by iteratively removing vertices of $G_{m,k}$ of degree at most $k/2$, one obtains a (non-empty) graph of minimum degree at least $k/2$.
Let $G$ be a connected component of the largest average degree in the obtained graph.
%Note that $\s{V(G)} \ge k/2$, and
Note that the average degree of $G$ is at least $k$, the maximum degree is still at most $6k$, and the girth is still at least $\log m/(6k)^2$.
Finally, $G$ still satisfies item $(iii)$ of Theorem~\ref{theo:cube}, being a subgraph of $G_{m,k}$.

Let $G_\circ$ be obtained by taking the replacement product of $G$ with a cycle. That is, $G_\circ$ is obtained from $G$ by replacing each vertex of degree $x$ by a cycle on $x$ new vertices -- which we henceforth refer to as a ``cloud'' -- and further adding edges as follows: if $u,v$ are adjacent in $G$, with $u$ being the $\ith$ neighbor of $v$ and $v$ being the $\jth$ neighbor of $u$ (under a fixed arbitrary enumeration of the neighbors of each vertex), then the $\ith$ vertex in the cloud corresponding to $v$ is connected by an edge to the $\jth$ vertex in the cloud corresponding to $u$.
So for example, it is easy to see that there is a one-to-one correspondence between the edges of $G$ and those edges of $G_\circ$ that connect vertices from different clouds.
Note that our graph $G_\circ$ is connected, as needed. Letting $n$ denote the number of its vertices,
note that $n$ equals the sum of the degrees of all vertices of $G$, so $n \ge k\s{V(G)} \ge n_0$, as needed.
Furthermore, $G_\circ$ is $3$-regular, since each vertex has two neighbors in its cloud and one neighbor in precisely one other cloud, as required by item~$(i)$ of the statement.

Let us now prove that the girth of $G_\circ$ is equal to the minimum between the girth of $G$ and the minimum degree of $G$. First, note that any cycle $C$ in $G_\circ$, other than a cloud, naturally determines a closed trail in $G$ (i.e., where vertices may be visited more than once, but not edges). %(since a pair of clouds have at most one edge between them
Indeed, for each edge of $C$ that connects two different clouds, the trail simply moves along the corresponding edge in $G$.. Note that the length of $C$ is at least the length (i.e., number of edges) of the trail. Since the length of the shortest closed trail in $G$ is its girth, we conclude that the length of any cycle in $G_\circ$ is at least the girth of $G$, unless that cycle is a cloud. Furthermore, since the smallest number of vertices in a cloud of $G_\circ$ equals the minimum degree of $G$, our claim follows.
That is, the girth of $G_\circ$ is at least
$$\min\{\log m/(6k)^2,\,k/2\} =\log\log(m)/1296 \ge \log\log(n)/2000 \;,$$
where we used the setting of $k$, Equation~(\ref{eq:relations}) and the fact that $n \le 6k\s{V(G)} \le 6km \le m^2$.
This proves item~$(ii)$ of the statement.

It remains to show that $G_\circ$ satisfies item~$(iii)$ of the statement as well.
Let $H_\circ$ be a $t$-vertex subgraph of $G_\circ$. Our goal is to bound $\phi_{H_\circ}$ from above.
Let $H$ be the induced subgraph of $G$ obtained by retaining only those vertices whose corresponding cloud has at least one vertex in $H_\circ$. Put $h=\s{V(H)}$, and notice $t \ge h$.
%We next apply item~(iii) of Theorem~\ref{theo:cube} which, as mentioned above, is also satisfied by $G$.
We next consider two cases, depending on whether $H$ is a forest or not.

First, suppose that $H$ is not a forest.
Hence, by item~$(iii)$ of Theorem~\ref{theo:cube} (a property which is also satisfied by $G$, as mentioned above) there is a partition $V(H)=S \cup S'$ with $\s{S},\s{S'} \ge h/3$ satisfying $\s{\partial_H(S)},\s{\partial_H(S')} \le (h/\log h) \cdot (\log\log h)^2$.
%Consider the subset of $V(H_\circ)$ corresponding to $S_1$ (i.e., obtained by replacing each vertex in $S_1$ with the vertices of its cloud in $H_\circ$) and to $S_2$.
Let $S_{\circ}$ be the subset of $V(H_\circ)$ corresponding to $S$ (i.e., obtained by replacing each vertex in $S$ with the vertices of its cloud in $H_\circ$). %Define $S'_{\circ}$ from $S'$ in a similar fashion.
Assume without loss of generality that $\s{S_{\circ}} \le t/2$ (otherwise take $S'_{\circ}$, which is defined from $S'$ in a similar fashion).
%Let $S \in \{S_1,S_2\}$ have the fewer vertices in its corresponding subset, and denote the corresponding subset by $S_{\circ} \sub V(H_\circ)$.
%Then $h/3 \le \s{S_{\circ}} \le t/2$.
%Observe that $\s{\partial_{H_\circ}(S_{\circ})} = \s{\partial_{H}(S)}$ since edges in $H$ precisely correspond to inter-cloud edges in $H_\circ$.
Observe that $\s{\partial_{H_\circ}(S_{\circ})} \le \s{\partial_{H}(S)}$, since any edge in $\partial_{H_\circ}(S_{\circ})$ must go between two different clouds, %(recall $S_{\circ}$ consists of whole clouds),
and there is a unique edge in $\partial_{H}(S)$ connecting the two vertices corresponding to these clouds.
Therefore,
$$\phi_{H_\circ} \le \frac{\s{\partial_{H_\circ}(S_{\circ})}}{\s{S_{\circ}}} \le
\frac{(h/\log h) \cdot (\log\log h)^2}{h/3} = \frac{3(\log\log h)^2}{\log h} \le \frac{3(\log\log t)^2}{\log (t/6k)} \le \frac{6(\log\log t)^2}{\log t} \;,$$
where %in the first inequality we used the fact that $\s{S_{\circ}} \le t/2$,
in the second inequality we used the fact that $\s{S_{\circ}} \ge \s{S} \ge h/3$, in the third inequality we used the fact that $h \le t \le 6k\cdot h$, and in the last inequality we used the fact that $t/6k \ge \sqrt{t}$ (i.e., $\sqrt{t} \ge 6k$); the latter follows from the fact that since $H$ is not a forest, $h$ is at least the girth of $G$, so $t \ge h \ge \log m/(6k)^2 \ge (6k)^2$ by Equation~(\ref{eq:relations}).
This proves item~$(iii)$ of the statement under the assumption that $H$ is not a forest.

Suppose next that $H$ is a forest. Notice we may assume that $H$ is a connected graph since otherwise $H_\circ$ is also not connected, meaning that $\phi_{H_\circ} = 0$ so we are done.
We apply Claim~\ref{clm:tree} on the tree $H$, where we set the weight of each vertex in $H$ to be the number of vertices in the corresponding cloud in $H_\circ$. Let $v$ be the vertex guaranteed by Claim~\ref{clm:tree}, and let $v_1,\ldots,v_d$ be the vertices of the cloud/cycle corresponding to $v$, in their order on the cycle. For each $1 \le i \le d$, let $S_i \sub V(H_\circ)$ be the set of vertices in $H_\circ$ corresponding to the $\ith$ connected components of $H-v$ (i.e., so that $v_i$ is the unique vertex in the cloud of $v$ that is connected to $S_i$).
Put $S'_i = S_i\cup\{v_i\}$. Then $\sum_{i=1}^d \s{S'_i} = t$, and our choice of $v$ guarantees that $\s{S'_i} \le t/2 + 1$.
We claim that there is an index $1 \le j \le d$ such that $(1/4)t \le \sum_{i=1}^j \s{S'_i} \le (3/4)t$.
Indeed, if $1\le j\le d$ is the smallest index such that $\sum_{i=1}^j \s{S'_i} \ge (1/4)t$ then
$$\sum_{i=1}^j \s{S'_i} = \sum_{i=1}^{j-1} \s{S'_i} + \lvert S'_j \rvert \le (t/4-1) + (t/2+1) = (3/4)t \;.$$
Now, let $S_\circ \sub V(H_\circ)$ be the smallest between $\bigcup_{i=1}^j S'_i$ and its complement, so that $t/4 \le \s{S_\circ} \le t/2$.
Observe that since $\{1,2,\ldots,j\}$ is an interval,  $\s{\partial_{H_\circ}(S_\circ)} \le 2$.
%As $t/4 \le \s{S_\circ} \le t/2$,
We conclude that
$$\phi_{H_\circ} \le 2/(t/4) = 8/t \le (1/\log t)\cdot (4\log\log t)^2 \;,$$
where in the last inequality we used the fact that $(x/\log x)\cdot (\log\log x)^2 \ge 1/2$, which can be verified to hold for any real $x\ge 3$ (and thus for any integer $t>2$).
This completes the proof.
%\end{proof}
\EPFOF

\subsection{Lower bound proof}\label{sec:LB-proof}

For our proof of Theorem~\ref{theo:main} we will need the graph witnessing the lower bound to contain a bridge.
The following lemma shows that one can modify a given graph so as to contain a bridge while preserving high girth and small hereditary expansion.

\BL\label{lemma:bridge} Suppose there is a $3$-regular
connected $n$-vertex graph $G$ with girth $g$ that is
$f$-non-expanding, where $f:[1/2,\infty)\to\R$ is monotone
decreasing. Then there is a $3$-regular connected
$(2n+2)$-vertex graph that contains a bridge, and moreover,
has girth at least $g$ and is $h$-non-expanding with
$h(x)=3f(x/2-1)$.
\EL
% \begin{proof}
\BPF Let $G_1,G_2$ be two vertex-disjoint copies of $G$.
Let $e_i$ be an arbitrary edge of $G_i$, $i\in\{1,2\}$, and
let $G_i'$ be obtained by subdividing $e_i$. That is,
$G_i'$ is obtained from $G_i$ by adding a new vertex $w_i$,
removing the edge $e_i=(u_i,v_i)$ and adding the edges
$(u_i,w_i),(w_i,v_i)$. It is clear that subdividing an edge
does not decrease the girth. Now, construct the graph $F$
from the union of $G_1'$ and $G_2'$ by adding the bridge
$(w_1,w_2)$. It is clear that $F$ is $3$-regular, connected
and has girth at least $g$. It therefore remains to show
that $F$ is $h$-non-expanding. Let $H$ be a $t$-vertex
subgraph of $F$ with $t>2$. We need to show that $\phi_H \le h(t)$.
Without loss of generality, $H$
has at least $t/2$ vertices in $G_1'$. Let $H'$ be the
subgraph of $H$ induced by those vertices, where we remove
the subdividing vertex $w_1$ if $w_1 \in V(H)$. Note that
$H'$ is a subgraph of $G_1$. Let $t' \ge t/2-1$ denote the
number of vertices of $H'$. Since $H'$ is
$f$-non-expanding, there is a subset $S \sub V(H')$ with
$\s{S} \le t'/2$ and $\s{\partial_{H'}(S)}/\s{S} \le f(t')$. Note
that $\s{\partial_{H}(S)} \le \s{\partial_{H'}(S)}+2$,
since the only edges in $H$ connecting a vertex in $H'$ and
a vertex not in $H'$ are $(u_1,w_1)$ and $(w_1,v_1)$. We
conclude that
$$\phi_H \le 3f(t') \le 3f(t/2 - 1) = h(t) \;,$$
where in the second inequality we used the monotonicity of $f$ for $t \ge 1/2$.
%This completes the proof.
%\end{proof}
\EPF

\medskip

%Henceforth, for a graph $G$ and a bijection $\sigma:V(G)\to[n]$, we write $\sigma(G)$ for the isomorphic graph to $G$ with vertex set $[n]$ and edges of the form $\sigma(e)=(\sigma(u),\sigma(v))$ for every $e=(u,v) \in E(G)$.
%
%$\sigma(G)=(\,[n],\,\{\sigma(e):e\in E(G)\}\,)$ be the input graph to $A$.\footnote{For an edge $e=(u,v)$ we write $\sigma(e)$ for the edge whose endpoints are $\sigma(u)$ and $\sigma(v)$.}

For a local sparse spanning graph algorithm $\calA$, we denote by $\calA(G,u,v) \in \{0,1\}$ the output of $\calA$ when the input graph is $G=(V,E)$ and the input edge is $(u,v) \in E$.
The \emph{query-answer transcript} of $\calA$ on $G$, where $\calA$ makes $q$ queries and $G$ is $d$-regular, is the sequence of triples $((x_j,i_j,y_j))_{j=1}^q$ where $(x_j,i_j) \in V\times[d]$ is the $\jth$ query and $y_j\in V$ is the corresponding answer.

%The \emph{query-answer transcript} of a local algorithm with query complexity $q$ is the sequence of triples $((x_j,i_j,y_j))_{j=1}^q$ where $(x_j,i_j) \in [n]\times[3]$ is the $\jth$ query and $y_j\in[n]$ is the corresponding answer.
Finally, for a permutation $\sigma$ on $V$,
%if $G$ is a graph whose vertex set is $[n]$ and $\sigma$ is a permutation on $[n]$ %then
we denote by $\sigma(G)$ the graph isomorphic to $G$ on the same vertex set, for which $(u,v) \in E(\sigma(G))$
if and only if $(\sigma(u),\sigma(v)) \in E(G)$.
We stress that in what follows, the graph $\sigma(G)$ will not necessarily have the same neighborhood ordering as that of $G$.
That is, if $y$ is the $\ith$ neighbor of $x$ in $G$ and $\sigma(v)=x,\sigma(u)=y$ then $u$ is {\em not} necessarily the $\ith$ neighbor of $v$ in $\sigma(G)$.

\BL
\label{lemma:SSG_girth}
Let $G$ be a $3$-regular
connected graph of girth $g$ that contains a bridge. Any
$(\frac12,q)$-local sparse spanning graph algorithm for the
graphs isomorphic to $G$ satisfies $q \ge g/2$.
%(for $\eps \leq 1/4$).
\EL
% \begin{proof}
\BPF
Let $\calA$ be an $(\frac12,q)$-local sparse
spanning graph algorithm for the graphs isomorphic to $G$,
and assume, contrary to the claim in the lemma, that $q < g/2$.
%$\calA$ has query complexity $q<g/2$ when run with $\eps < 1/2$.
%$\eps = 1/4$.
We shall say that $\calA$ {\em accepts\/}
an edge $(u,v)$ in $G$ if it gives a positive answer when queried on $(u,v)$ (that is, $(u,v)$ belongs to the sparse spanning graph $G'$).
We will show that with probability $1$
over its random coins, $\calA$ accepts every edge of $G$. This
will complete the proof as it  means that the number
of edges of $G$ that $\calA$ accepts is $(1+\frac12)n$,
%$3n/2 > (1+\frac14)n$,
where $n$ is the number of vertices of $G$,
contradicting condition~$(ii)$ in
Definition~\ref{dfn:SSG-alg}.

Let $(u,v) \in E(G)$ and assume for contradiction that there is a sequence $r$ of random coins for $\calA$ such that the
corresponding deterministic algorithm $\calA_r$ satisfies $\calA_r(G,u,v)=0$.
Suppose, without loss of generality, that the vertex
set of $G$ is $[n]$ and that
$(1,2)$ is a bridge in $G$.
We will construct a permutation $\sigma$ on $[n]$ with $\sigma(u)=1,\, \sigma(v)=2$
so that the graph $\sigma(G)$ (with an appropriate way of ordering the neighbors of each vertex)
has the property that the query-answer transcript of $\calA_r(G,u,v)$ is identical to that of $\calA_r(\sigma(G),u,v)$.
Note that $\calA_r(\sigma(G),u,v)$ is well defined since the input edge $(u,v)$ is indeed an edge of $\sigma(G)$, and since $\sigma(G)$ is a valid input graph to $\calA_r$ being isomorphic to $G$.
Since $\calA_r$ is deterministic, whether or not $\calA_r$ accepts $(u,v)$ depends solely on the
query-answer transcript.
Therefore, the existence of $\sigma$ as above would imply that $\calA_r(\sigma(G),u,v)=0$. However, this would contradict condition~$(i)$ in Definition~\ref{dfn:SSG-alg} since $(u,v)$ is a bridge in $\sigma(G)$.

Let $Q=(x_j,i_j,y_j)_{j=1}^q$ be the query-answer transcript of $\calA_r(G,u,v)$.
We first claim that if a permutation $\sigma$ and an ordering of the neighbors of each vertex
of $\sigma(G)$, are such that $\sigma(u)=1,\sigma(v)=2$ and
for every $1 \leq j \leq q$ the $i_j$-th neighbor of vertex $x_j$ in $\sigma(G)$ is vertex $y_j$ then the query-answer transcript of $\calA_r(G,u,v)$ is
identical to the query-answer transcript of $\calA_r(\sigma(G),u,v)$.
To see this, let the query-answer transcript of $\calA_r(\sigma(G),u,v)$ be denoted by $(x_j',i_j',y_j')_{j=1}^{q'}$.
We prove, by induction on $j$, that the two query-answer transcripts are the same when restricted to the first $1\le j\le q$ queries, that is, $(x'_j,i'_j)=(x_j,i_j)$ and $y'_j=y_j$ for every $1\le j\le q$.
Note that this will also imply that $q=q'$ (i.e., that the number of queries is identical).
For $j=1$ we have $(x'_1,i'_1)=(x_1,i_1)$ since $\calA_r$ is deterministic and in both cases the input is $(u,v)$.
Our assumption on $\sigma$ thus guarantees that we also have $y'_1=y_1$.
Suppose our claim holds for the first $j-1$ queries.
Again, since $\calA_r$ is deterministic, the $\jth$ query is determined only by the query-answer transcript of the first $j-1$ queries (and the input edge).
Hence, the induction hypothesis implies that $(x'_j,i'_j)=(x_j,i_j)$ and our assumption on $\sigma$
again implies that we also have $y'_j=y_j$. This completes the inductive proof.

%\mnote{D: modified the remainder of the proof a bit}
It follows that in order to
complete the proof it suffices to find a permutation $\sigma$ and an ordering of
the neighbors of each vertex, as above.
Let again $Q=(x_j,i_j,y_j)_{j=1}^q$ be the query-answer transcript of $\calA_r(G,u,v)$,
and let $F$ be the (labeled) graph spanned by the edge set\footnote{$E(F)$ might contain the edge $(x,y)$ twice if $y$ is the $\ith$ neighbor of $x$, $x$ is the $\jth$ neighbor of $y$ and the algorithm queried both $(x,i)$ and $(y,j)$. In this case we will keep just one copy of $(x,y)$ thus  making sure
that $E(F)$ is indeed a set, and not a multi-set.}
$$E(F)= \{(x_j,y_j)\}_{j=1}^q \cup \{(u,v)\} \;.$$
Since
\begin{equation}\label{eq:EF}
\s{E(F)} \le q+1 \le g/2 \;,
\end{equation}
we have that $F$ is a forest.
Let $T_1,\dots,T_k$ be the (labeled) trees in $F$. For the sake of defining
$\sigma$ it will be convenient to consider a single tree $T$.
% If $F$ is composed of $k$ labeled trees $T_1,\ldots,T_k$
The edge-set of $T$ consists of $E(F)$ and $k-1$ additional
edges. The additional edges do not necessarily belong to $G$, and are selected as follows.
For each labeled tree $T_i$, let $t_i$ denote an
arbitrary vertex of degree smaller than $3$. For every
$i\in [k-1]$, add the edge $(t_i, t_{i+1})$.

Observe that Equation~(\ref{eq:EF}) implies that $\s{E(T)} <
g$. Consider a rooted version of $T$ where $u$ is the root,
and construct $\sigma$ as follows. Set $\sigma(u)=1$,
$\sigma(v)=2$, and define the neighborhood relation between
$u,v$ in $\sigma(G)$ as it is in $G$. That is, if $u$ is
the $\ith$ neighbor of $v$ and $v$ is the $\jth$ neighbor
of $u$ in $G$ then the same holds in $\sigma(G)$. Suppose
we have already defined $\sigma(x)$ for all $x$ at distance
at most $d-1$ from $u$ (in $T$) as well as for some
vertices at distance $d$, and let $y$ be a vertex at
distance $d$ for which $\sigma(y)$ has not been defined
yet. Let $x$ be the parent of $y$ in $T$ (whose distance
from $v$ is thus $d-1$) and let us set $\sigma(y)$ to be a
neighbor of $\sigma(x)$ in $G$ which is not the image of
any vertex under the $\sigma$ we have defined thus far.
Such a vertex exists since $G$ is $3$-regular and the
degree in $T$ is at most $3$. If the edge $(x, y)$ is in
$F$ then we define the neighborhood relation between
$\sigma(x)$ and $\sigma(y)$ as $x$ and $y$ in $G$. Once we
define $\sigma$ for all vertices of $T$ we arbitrarily
extend $\sigma$ to a permutation, and extend the
neighborhood relation between the vertices in a consistent
manner.
% It is also easy to see that we can extend the
% neighborhood relation between the vertices in some way.
\EPF

We are now ready to prove Theorem~\ref{theo:main}.

\BPFOF{Theorem~\ref{theo:main}}
Let
$$h(x)=\frac{1}{\log(3x)}\cdot 32(\log\log(8x))^2 \;.$$
It is not hard to check that $h:[1/2,\infty)\to\R$ is monotone decreasing.
Note that the graph in Lemma~\ref{lemma:cube-cycle} is $h$-non-expanding, since for $x\ge 3$,
$$h(x) \ge \frac{1}{\log(x^2)}\cdot 32(\log\log x)^2 = \frac{1}{\log x}\cdot (4\log\log x)^2 \;.$$
Apply Lemma~\ref{lemma:bridge} on the graph(s) in Lemma~\ref{lemma:cube-cycle}.
We get a $3$-regular connected $n$-vertex graph, for infinitely many $n$, that contains a bridge, has girth at least
$$\frac{\log\log((n-1)/2)}{2000} \ge \frac{\log\log(n/4)}{2000} \ge \frac{\log\log(n)}{4000}$$
and is $f$-non-expanding with
$$f(x)= 3h(x/2-1) \le \frac{1}{\log (x/2)} \cdot 96(\log\log(4x))^2
\le \frac{3}{\log x} \cdot 96(4\log\log x)^2
\le \frac{1}{\log x} \cdot (70\log\log x)^2$$
where we assumed $x\ge 3$.
The proof now follows immediately from Lemma~\ref{lemma:SSG_girth}.
\EPFOF

\bibliographystyle{plain}
%\bibliography{refs}

\end{document}